\documentclass[a4paper,draft]{amsproc}
\usepackage{amssymb}
\usepackage[hyphens]{url} 

\theoremstyle{plain}
 \newtheorem{theorem}{Theorem}[section]

 \newtheorem{corollary}{Corollary}[section]
\theoremstyle{definition}

\theoremstyle{remark}

 \numberwithin{equation}{section}

\renewcommand{\leq}{\leqslant}
\renewcommand{\geq}{\geqslant}

\title[partition polynomials / M. Goubi]{On partition polynomials and partition functions}



\author{\bfseries Mouloud  Goubi} 
\address{Mouloud Goubi\\
Department of Mathematics \\
University of UMMTO RP. 15000\\
Tizi-ouzou, Algeria\\
Laboratoire d'Alg\`ebre et Th\'eorie des Nombres, USTHB Alger}
\email{mouloud.goubi@ummto.dz}

\begin{document}

\vspace{18mm} \setcounter{page}{1} \thispagestyle{empty}

\begin{abstract}
In this paper we revisit the work of E.T. Bell concerning partition
polynomials in order to introduce the reciprocal partition
 polynomials. We give their explicit formulas and apply the result to
 compute closed formulae for some well-known partition
 functions.

\end{abstract}

\maketitle

\section{Introduction}
Many known and new arithmetical functions are included as special
cases of partition polynomials. In this work we consider partition
polynomials $P_n(z)$ introduced and studied by E. T. Bell in the
work \cite{BE}. We introduce and study the reciprocal polynomials
$W_n(z)$ which are a generalization of partition function $p(n)$ and
restricted partition functions $W\left(n,d^s\right).$ Formal
calculus allow us to compute $W_n(z)$ and deduce explicit formula
for $W\left(n,d^s\right)$ and $p(n)$. We end the work by the
generalized partition polynomials $WP_n(z)$ including $P_n(n)$ and
$W_n(z)$. The result conducts to explicit formula of a large family
of partition functions.

\section{Partition polynomials and properties}
We reproduce here the family of polynomials constructed by E.T. Bell
in the work \cite{BE}. Let $n>0$ be an integer and $C_j$ denote a
set of distinct integers $>0$. $C_j$ may contain any finite or
infinite number of elements. The polynomial $\psi^{(j)}_n(z_j)$ for
a given $z_j$ and $C_j$ is defined by
\begin{equation}
\psi^{(j)}_n(z_j)=\sum_{d\in C_j\atop d|n}dz^{n/d}_j.
\end{equation}
We consider now $j=1,2,\cdots,s$ and $a_j$ are integers not all
zero. Let $z=\left(z_1,\cdots,z_s\right)$, the polynomial
$\Psi_n\left(z\right)$ is defined by
\begin{equation}
\Psi_n\left(z\right)=-\sum_{j=1}^{s}a_j\psi^{(j)}_n(z_j).
\end{equation}
The partition polynomial $P_n(z,C,a)=P_n(z)$ of rank $n$, argument
$z$, set $C=\left\{C_1,C_2,\cdots,C_s\right\}$ and index
$a=\left(a_1,a_2,\cdots,a_s\right)$ is given by
\begin{equation}\label{pnz}
P_n(z)=\sum_{\pi(n)}\prod_{j=1}^{n}\left[\frac{1}{k_j!}\left(\frac{\Psi_j(z)}{j}\right)^{k_j}\right],
\end{equation}
where
$\pi(n)=\left\{\left(k_1,\cdots,k_n\right)\in\mathbb{N}\backslash
k_1+2k_2+\cdots+nk_n=n\right\}.$ Some recursive formulae for
$P_n(z)$ are
\begin{equation}
P_n(z,C,a+b)=\sum_{j=0}^{n}P_j(z,C,a)P_{n-j}(z,C,b),
\end{equation}
\begin{equation}
P_n(z,A+B,a+b)=\sum_{j=0}^{n}P_j(z,A,a)P_{n-j}(z,B,b).
\end{equation}
We say that $f(t)$ is a generating function for the sequence
$\left(c_n\right)_{n\in\mathbb{N}}$ of numbers or polynomials if
$f(t)$ is written as a power series, then
\begin{equation}\label{pser}
f(t)=\sum_{n\geq0}c_nt^n.
\end{equation}
We do not have to concern ourselves with questions of convergence of
the series \eqref{pser}, since we are interested in the coefficients
$c_n$. We consider such series as formal power series in $t$; for
more information about this theory we refer to  account given by
Ivan Niven \cite{NI}. In the same sense the generating function of
polynomials $P_n(z)$ is
\begin{equation}
f(t)=f\left(t,z,C,a\right)=\prod_{1}\left(1-z_1t^{n_1}\right)^{a_1}\prod_{2}\left(1-z_2t^{n_2}\right)^{a_2}\cdots\prod_{s}\left(1-z_st^{n_s}\right)^{a_s};
\end{equation}
where $\prod_j$ denotes a product with respect to all $n_j$ such
that $n_j\in S_j$. E.T. Bell gave the croquet of the proof by
considering $f(t)=\exp\log f(t)$ and indicates to use Maclaurin's
theorem to get corresponding series expansion. Here we revisit the
proof by using Fa\`{a} di Bruno formula (see \cite{FA}). If $h(t)$
and $g(t)$ are functions for which all the necessary derivatives are
defined, then
\begin{equation}
\left(h\circ
g\right)^{(n)}(t)=\displaystyle\sum_{k=k_1+\cdots+k_n\atop
k_1+2k_2+\cdots+nk_n=n}\frac{n!}{k_1!\cdots
k_n!}h^{(k)}(g(t))\prod_{i=1}^{n}\left(\frac{g^{(i)}(t)}{i!}\right)^{k_i};
\end{equation}
where $f^{(n)}(t)=\frac{d^nt}{dt^n}$. A detailed proof is given by
Steven Roman \cite{RO} by using the umbral calculus. Then the
coefficients $[t^n]h\circ g(t)$ of the series expansion of $h\circ
g(t)$ take the form $[t^0]h\circ g(t)=h(g(0))$ and for $n\geq1$;
\begin{eqnarray}
[t^n]h\circ g(t)=\displaystyle\sum_{k=k_1+\cdots+k_n\atop
k_1+2k_2+\cdots+nk_n=n}\frac{h^{(k)}(g(0))}{k_1!\cdots
k_n!}\prod_{i=1}^{n}\left(\frac{g^{(i)}(0)}{i!}\right)^{k_i}.
\end{eqnarray}

We know that
\[f(t)=\exp\left(\sum_{j=1}^{s}\sum_{n_j\in
S_j}a_j\log\left(1-z_jt^{n_j}\right)\right),\] but we have
$\log\left(1-t\right)=-\sum_{n\geq1}\frac{t^k}{k}.$ Then
\[f(t)=\exp\left(-\sum_{k\geq1}\frac{1}{k}\sum_{j=1}^{s}\sum_{n_j\in S_j}a_jz^{k}_jt^{n_jk}\right).\]
For computing the formal power series of $f$ let
\[g(t)=-\sum_{k\geq1}\frac{1}{k}\sum_{j=1}^{s}\sum_{n_j\in
S_j}a_jz^{k}_jt^{n_jk},\] then $f(t)=e^{g(t)}.$ Furthermore for
$n\geq1$ we have
\[[t^n]f(t)=\sum_{k=1}^{n}\sum_{k_1+\cdots+k_n=k\atop 1k_1+2k_2+\cdots+nk_n}\frac{n!}{k_1!\cdots k_n!}\prod_{i=1}^{n}
\left(\frac{-\sum_{j=1}^{s}\sum_{n_j\in S_j\atop
n_j|i}a_jn_jz^{i/n_j}_j}{i}\right)^{k_i}.\] Since
\[\sum_{j=1}^{s}\sum_{n_j\in S_j\atop
n_j|i}a_jn_jz^{i/n_j}_j=\sum_{j=1}^{s}a_j\psi^{(j)}_i(z_j)=\Psi_i(z),\]
then
\begin{equation}
P_n(z)=\sum_{k=1}^{n}\sum_{k_1+\cdots+k_n=k\atop
k_1+2k_2+\cdots+nk_n=n}\prod_{j=1}^{n}\left[\frac{1}{k_j}\left(\frac{\Psi_j(z)}{j}\right)^{k_j}\right],
\end{equation}
identic to the expression \eqref{pnz} bellow.

\section{Explicit formula of reciprocal partition polynomials}
Let us introducing the family $W_n(z)=W_n(z,C,a)$ of partition
polynomials, where $W_n(z)$ is generated by the function $1/f(t).$
The Cauchy product of generating functions of $P_n(z)$ and $W_n(z)$
equal $1$. Hence $W_0(z)=1$ and others are obtained from the
recursive formula
\begin{equation}
W_n(z)=-\sum_{k=0}^{n-1}W_k(z)P_{n-k}(z).
\end{equation}
For the prove we refer to \cite{Gou}. From the definition of
$W_n(z)$, we can write $W_n(z)=P_n(z,C,-a)$ with
$-a=\left(-a_1,\cdots,-a_s\right)$. According to the relation
$1/f(t)=\exp\left(-\log f(t)\right)$, the following theorem is
immediate
\begin{theorem}
\begin{equation}\label{wnz}
W_n(z)=\sum_{\pi(n)}(-1)^{\sum
k_j}\prod_{j=1}^{n}\left[\frac{1}{k_j!}\left(\frac{\Psi_j(z)}{j}\right)^{k_j}\right].
\end{equation}
\end{theorem}
We can prove the identity \eqref{wnz} with another method; which is
based on exponential partial Bell polynomials. To learn more about
this technique we refer to recent works \cite{Go,Goub,Goubi}.

\subsection{Expression of restricted partition function $W(n,d^s)$}
Special case, namely, restricted partition function
\[W(n,d^s)=W\left(n,\left\{d_1,d_2,\cdots,d_s\right\}\right)\] is completely studied, but the given formulas still so much big.
$W(n,d^s)$ is a number of partitions of $n$ into positive integers
$d_1,d_2,\cdots,d_s$ each not greater than $s$. The corresponding
generating function has the form
\begin{equation}
\prod_{i=1}^{s}\frac{1}{1-x^{d_i}}=\sum_{n\geq0}W(n,d^s)x^n.
\end{equation}
$W(n,d^s)$ satisfies the basic recursive relation
\begin{equation}
W(n,d^s)-W(n-d_s,d^s)=W(n,d^{s-1}).
\end{equation}
Sylvester (\cite{SY, SYL})showed that the restricted partition
function may be presented as a sum of Sylvester waves
\begin{equation}
W(n,d^s)=\sum_{j}W_j(n,d^s),
\end{equation}
where the sum $\sum_{j}$ is over all distinct factors of the
elements in the set $d^s$. B.Y. Rubinstein and L.G. Fel (see
\cite{RU}) proved that
\begin{eqnarray}
W_j(n,d^s)&=&\frac{1}{(\omega_j-1)!\pi_{\omega_j}}\sum_{\rho_j}\frac{\rho^{-n}_j}{\prod_{i=\omega_j+1}\left(1-\rho^{d_i}_j\right)}\times\\
\nonumber&&\sum_{k=0}^{\omega_j-1}B^{(w_j)}_k\left(n+n_{w_j}|d^{w_j}\right)H^{(s-w_j)}_{w_j-1-k}\left(n_s-n_{w_j},\rho_{j}|d^{s-\omega_j}\right),
\end{eqnarray}
where
\[\frac{e^{st}\prod_{i=1}^{m}\left(1-\rho^{d_i}\right)}{\prod_{i=1}^{m}e^{d_it}-\rho^{d_i}}=\sum_{n\geq0}H^{(m)}_{n}\left(\rho\mid d_n\right)\frac{t^n}{n!},\
\rho^{d_i}\neq1.\] and
\[\frac{e^{st}t^m\prod_{i=1}^{m}d^i}{\prod_{i=1}^{m}\left(e^{d_it}-1\right)}=\sum_{n\geq0}B^{(m)}_{n}\left(s|d^m\right)\frac{t^n}{n!}.\]
$W(n,d^s)$ corresponds to $W_n(z)$ in the case $z=a={\bold
1}=\left(1,1,...,1\right)$ and $C_i=\left\{d_i\right\}$ for all
$i=(1,2,\cdots,s)$. According to these conditions we have for all
$j=1,\cdots,s$:
\[\Psi_n\left({\bold1}\right)=-\sum_{1\leq j\leq s\atop d_j|n}d_j.\]
We define the restricted divisor function to $S$; $d_{S}$ by
$d_{S}(j)=\sum_{d_i|j}d_i$, then a simple formula of
$W\left(n,d^s\right)$ is given by the following theorem.
\begin{theorem}
\begin{equation}
W\left(n,d^s\right)=\sum_{\pi(n)}\prod_{j=1}^{n}\left[\frac{1}{k_j!}\left(\frac{d_S(j)}{j}\right)^{k_j}\right].
\end{equation}
\end{theorem}

\subsection{Generating functions of partitions}
Let a set $S\subset\mathbb{N}$ and $p(n|S)$ the number of partitions
of $n$ into elements of $S$. Then the generating function of
$p(n|S)$ is
\begin{equation}
\prod_{k\in S}\frac{1}{1-t^k}=\sum_{n\geq0}p(n|S)t^n.
\end{equation}
If $p_m(n|S)$ is the number of partitions with exactly m-part in
$S$, then the generating function is
\begin{equation}
\prod_{k\in S}\frac{1}{1-xt^k}=\sum_{m,n\geq0}p_m(n|S)x^mt^n.
\end{equation}
In fact we have
\[\prod_{k\in S}\frac{1}{1-t^k}=\prod_{k\in S}\sum_{m_i\geq0}t^{km}=\sum_{m\geq0}t^{\sum_{k\in S}mk}.\]
Then \[[t^n]\displaystyle\prod_{k\in
S}\frac{1}{1-t^k}=\sum_{\displaystyle\sum_{k\in S\atop
m\geq0}mk=n}1.\] For the second, we have
\[\prod_{k\in S}\frac{1}{1-xt^k}=\prod_{k\in S}\sum_{m\geq0}x^mt^{mk}=\displaystyle\sum_{m_k\geq0\atop k\in S}x^{\sum_{k\in S}m_k}t^{\sum_{k\in S}km_k}.\]
Then \[[x^mt^n]\prod_{k\in S}\frac{1}{1-xt^k}=\sum_{\sum_{k\in
S}m_k=m\atop \sum_{k\in S}km_k=n}1.\] When $S=\mathbb{N}$, the
corresponding generating functions may be displayed, respectively,
as
\begin{equation}\label{pn}
\prod_{k=1}^{\infty}\frac{1}{1-t^k}=\sum_{n\geq0}p(n)t^n
\end{equation}
and
\begin{equation}
\prod_{k=1}^{\infty}\frac{1}{1-xt^k}=\sum_{m,n\geq0}p_m(n)x^mt^n.
\end{equation}
The elementary aspects of the theory of partitions are given in
detail in \cite[Chap.19]{HA}. A partition of a positive integer $n$
may be thought as an unordered representation of $n$ as a sum of
other positive integers. Thus $3+2, 2+3$ represent the same
partition of $5$. Euler gave the first recurrent formula
(see\cite{EU}) of the arithmetical function $p(n)$:
\begin{equation}
p(n)=\sum_{k=1}^{n}(-1)^{k+1}\left[p\left(n-\frac{1}{2}k\left(3k-1\right)\right)+p\left(n-\frac{1}{2}\left(3k+1\right)\right)\right].
\end{equation}
Andrews after proving this formula (see \cite{AN}), he says "No one
has ever found a more efficient algorithm for computing $p(n)$. It
computes a full table of values of $p(n)$ for $n > 5$, in time
$O(n^{3/2} )$." In 1917 Hardy and Ramanujan (see \cite{HAR}) applied
on \eqref{pn} the theory of functions of complex variables and
developed a method which yields an asymptotic formula for $p(n)$:
\begin{eqnarray}\label{rama}
p(n)=\frac{1}{2\pi\sqrt{2}}\sum_{k\leq\alpha\sqrt{n}}A_k(n)\frac{d}{dn}\left(\frac{\exp\left(\frac{C\sqrt{n-1/24}}{k}\right)}{\sqrt{n-1/24}}\right)+O(n^{-1/4}),
\end{eqnarray}
with $\alpha$ as an arbitrary constant, and
\begin{equation}
A_k(n)=\sum_{h\mod k\atop (h,k)=1}\omega_{h,k}e^{-2\pi ihn/k},
C=\pi\sqrt{2/3}.
\end{equation}
Rademacher (see \cite{RA}) replaced the asymptotic formula
\eqref{rama} by the equality
\begin{eqnarray}
p(n)=\frac{1}{2\pi\sqrt{2}}\sum_{k\geq1}A_k(n)\frac{d}{dn}\left(\frac{\exp\left(\frac{C\sqrt{n-1/24}}{k}\right)}{\sqrt{n-1/24}}\right),
\end{eqnarray}
in which the series is absolutely convergent. Recently; Aleksa
Srdanov (see \cite{SR}) investigated the arithmetical function
$p(n)$. First he defines numbers $p(n,m)$ the number of all possible
partitions of the number $n$ having exactly $m$ parts, ($1\leq m\leq
n$). Then $p(n)=\sum_{k=1}^{n}p(n,k)$. Finally he computed in
different way the expression of numbers $p(n,k)$; for more details,
we refer to Theorems 1,2,3 and 4 in the
work \cite{SR}.\\
For $S=\mathbb{N}$, we have \[d_S(j)=\sigma(j)=\sum_{i|j}i\] and the
following theorem is immediate
\begin{theorem}
We have $p(0)=1$ and for $n\geq1$;
\begin{equation}
p(n)=\sum_{\pi(n)}\prod_{j=1}^{n}\left[\frac{1}{k_j!}\left(\frac{\sigma(j)}{j}\right)^{k_j}\right].
\end{equation}
\end{theorem}
If there is some difficulties for computing $\sigma(n)$, we purpose
the following formula:
\begin{equation}
\sigma(n)=\prod_{j=1}^{m}\left[\sum_{k=0}^{\left\lfloor\frac{b_j}{2}\right\rfloor}(-1)^k{b_j-k\choose
k}p^{k}_j\left(1+p_j\right)^{b_{j}-2k}\right],
\end{equation}
when we know the decomposition of $n$ on prime factors;
$n=p^{b_1}_1\cdots p^{b_m}_m.$

\section{Generalized partition polynomials}
Let $w=\left(w_1,\cdots,w_l \right)\in\mathbb{C}^{l}$,
$b=\left(b_1,\cdots,b_l\right)\in\mathbb{N^{\star}}^{l}$ and
$S=\left(S_1,\cdots S_l\right)$; where $S_j$ contains finite or
infinite number of elements. We introduce the generalized partition
polynomials $WP_n(w,z)$ including polynomials $P_n(z)$ and $W_n(z)$
by the generating function:
\begin{equation}
h(t)=\frac{f\left(t,z,C,a\right)}{f\left(t,w,S,b\right)}.
\end{equation}
The following corollary is immediate
\begin{corollary}
\begin{equation}
WP_n(w,z)=\sum_{m=0}^{n}P_m(z)W_{n-m}(w).
\end{equation}
\end{corollary}

\subsection{Application to some partition functions}
In the literature finitely many partition functions are studied. We
focus our interest in partition functions $a(n)$, $\bar{a}(n)$,
$\psi^{\star}(n)$ and $\varphi^{\star}(n)$. The arithmetical
function $a(n)$ counts the number of partitions of weight $n$ such
that the even parts can appear in two colors (see \cite{CH,SE}). So,
for example, $a(3) = 4$ where the colored partitions in question are
\[3, 2_1+1, 2_2+1\ \textrm{and}\ 1+1+1.\] Its generating function
takes the form
\begin{equation}
\frac{1}{\prod_{n=1}^{\infty}\left(1-t^n\right)\prod_{n=1}^{\infty}\left(1-t^{2n}\right)}=\sum_{n\geq0}a(n)t^n.
\end{equation}
Chan (see \cite{CH}) proved that
\begin{equation}
3\frac{\prod_{n=1}^{\infty}\left(1-t^{3n}\right)^3\prod_{n=1}^{\infty}\left(1-t^{6n}\right)^3}{\prod_{n=1}^{\infty}\left(1-t^n\right)^{4}
\prod_{n=1}^{\infty}\left(1-t^{2n}\right)^{4}}
=\sum_{n\geq0}a\left(3n+2\right)t^n.
\end{equation}
Byungchan Kim (see \cite{KI}) introduced the overcubic partition
function $\bar{a}(n)$, which counts all of the overlined versions of
the cubic partitions counted by $a(n)$. Its generating function
takes the form
\begin{equation}
\frac{\prod_{n=1}^{\infty}\left(1-t^{4n}\right)}{\prod_{n=1}^{\infty}\left(1-t^{n}\right)^2\prod_{n=1}^{\infty}\left(1-t^{2n}\right)}=\sum_{n\geq0}\bar{a}(n)t^n.
\end{equation}
Kim provided that
\begin{equation}
6\frac{\prod_{n=1}^{\infty}\left(1-t^{3n}\right)^6\prod_{n=1}^{\infty}\left(1-t^{4n}\right)^3}{\prod_{n=1}^{\infty}\left(1-t^n\right)^{8}
\prod_{n=1}^{\infty}\left(1-t^{2n}\right)^{3}}
=\sum_{n\geq0}\bar{a}\left(3n+2\right)t^n.
\end{equation}
The Ramanujan's $\psi$ and $\varphi$ functions are defined as
\begin{equation}
\psi(t):=\sum_{n\geq0}t^{n(n+1)/2}
\end{equation}
and
\begin{equation}
\varphi(t)=1+2\sum_{n\geq1}t^{n^2}.
\end{equation}
These functions admit the following reformulations
\begin{equation}
\psi(t)=\frac{\prod_{n=1}^{\infty}\left(1-t^{2n}\right)^2}{\prod_{n=1}^{\infty}\left(1-t^n\right)}
\end{equation}
and
\begin{equation}
\varphi(t)=\frac{\prod_{n=1}^{\infty}\left(1-t^{2n}\right)^5}{\prod_{n=1}^{\infty}\left(1-t^{n}\right)^2\prod_{n=1}^{\infty}\left(1-t^{4n}\right)^2}
\end{equation}
$\psi^{\star}(n)$ and $\varphi^{\star}(n)$ the partition functions
generated respectively by $\psi$ and $\varphi$. The generating
functions of these partition functions are special case of the
function
\begin{eqnarray}
F(t)=\frac{\prod_{n=1}^{\infty}\left(1-t^{r_1n}\right)^{a_1}\prod_{n=1}^{\infty}\left(1-t^{r_2n}\right)^{a_2}}{\prod_{n=1}^{\infty}\left(1-t^{s_1n}\right)^{b_1}
\prod_{n=1}^{\infty}\left(1-t^{s_2n}\right)^{b_2}}
\end{eqnarray}
Let the function $I_i$ such that $I_i(j)=1$ if $i|j$ and zero
otherwise. We consider $WP(n)$ the partition function generated by
the function $F(t)$. According to Theorem we conclude that
\begin{eqnarray}
\nonumber
WP(n)&=&\sum_{\pi(n)}\prod_{j=1}^{n}\left[\frac{1}{k_j!}\left(\frac{b_1I_{s_1}(j)\sigma(j/s_1)+b_2I_{s_2}(j)\sigma(j/s_2)}{j}\right)^{k_j}\right]\\
\nonumber&+&\sum_{m=1}^{n}\sum_{\pi(m)}\prod_{j=1}^{m}\left[\frac{1}{k_j!}\left(\frac{-a_1I_{r_1}(j)\sigma(j/r_1)-a_2I_{r_2}(j)\sigma(j/r_2)}{j}\right)^{k_j}\right]\\
&\times&\sum_{\pi(n-m)}\prod_{j=1}^{n-m}\left[\frac{1}{k_j!}\left(\frac{b_1I_{s_1}(j)\sigma(j/s_1)+b_2I_{s_2}(j)\sigma(j/s_2)}{j}\right)^{k_j}\right]
\end{eqnarray}
From this identity follow the explicit formula for considered
partition functions:
\begin{equation}
a(n)=\sum_{\pi(n)}(-1)^{\sum
k_j}\prod_{j=1}^{n}\left[\frac{1}{k_j!}\left(\frac{\sigma(j)+I_2(j)\sigma\left(j/2\right)}{j}\right)^{k_j}\right],
\end{equation}
\begin{eqnarray}
\nonumber a(3n+2)&=&3\sum_{\pi(n)}4^{\sum
k_j}\prod_{j=1}^{n}\left[\frac{1}{k_j!}\left(\frac{\sigma(j)+I_2(j)\sigma\left(j/2\right)}{j}\right)^{k_j}\right]\\
\nonumber&+&3\sum_{m=1}^{n}\sum_{\pi(m)}(-3)^{\sum
k_j}\prod_{j=1}^{m}\left[\frac{1}{k_j!}\left(\frac{I_3(j)\sigma(j/3)+I_6(j)\sigma(j/6)}{j}\right)^{k_j}\right]\\
&\times&\sum_{\pi(n-m)}4^{\sum
k_j}\prod_{j=1}^{n-m}\left[\frac{1}{k_j!}\left(\frac{\sigma(j)+I_2(j)\sigma\left(j/2\right)}{j}\right)^{k_j}\right],
\end{eqnarray}
\begin{eqnarray}
\nonumber\bar{a}(n)&=&\sum_{\pi(n)}\prod_{j=1}^{n}\left[\frac{1}{k_j!}\left(\frac{2\sigma(j)+I_2(j)\sigma(j/2)}{j}\right)^{k_j}\right]\\
\nonumber&+&\sum_{m=1}^{n}\sum_{\pi(m)}(-1)^{\sum
k_j}\prod_{j=1}^{m}\left[\frac{1}{k_j!}\left(\frac{I_4(j)\sigma(j/4)}{j}\right)^{k_j}\right]\\
&\times&\sum_{\pi(n-m)}\prod_{j=1}^{n-m}\left[\frac{1}{k_j!}\left(\frac{2\sigma(j)+I_2(j)\sigma(j/2)}{j}\right)^{k_j}\right],
\end{eqnarray}
\begin{eqnarray}
\bar{a}\left(3n+2\right)&=&6\sum_{\pi(n)}\prod_{j=1}^{n}\left[\frac{1}{k_j!}\left(\frac{8\sigma(j)+3I_{2}(j)\sigma(j/2)}{j}\right)^{k_j}\right]\\
\nonumber&+&\sum_{m=1}^{n}\sum_{\pi(m)}\prod_{j=1}^{m}\left[\frac{1}{k_j!}\left(\frac{-6I_{3}(j)\sigma(j/3)-3I_{4}(j)\sigma(j/4)}{j}\right)^{k_j}\right]\\
&\times&\sum_{\pi(n-m)}\prod_{j=1}^{n-m}\left[\frac{1}{k_j!}\left(\frac{8\sigma(j)+3I_{2}(j)\sigma(j/2)}{j}\right)^{k_j}\right],
\end{eqnarray}
\begin{eqnarray}
\nonumber\psi^{\star}(n)&=&\sum_{\pi(n)}\prod_{j=1}^{n}\left[\frac{1}{k_j!}\left(\frac{\sigma(j)}{j}\right)^{k_j}\right]\\
\nonumber&+&\sum_{m=1}^{n}\sum_{\pi(m)}(-2)^{\sum
k_j}\prod_{j=1}^{m}\left[\frac{1}{k_j!}\left(\frac{I_2(j)\sigma(j/2)}{j}\right)^{k_j}\right]\\
&\times&\sum_{\pi(n-m)}\prod_{j=1}^{n-m}\left[\frac{1}{k_j!}\left(\frac{\sigma(j)}{j}\right)^{k_j}\right],
\end{eqnarray}
\begin{eqnarray}
\nonumber\varphi^{\star}\left(n\right)&=&\sum_{\pi(n)}2^{\sum
k_j}\prod_{j=1}^{n}\left[\frac{1}{k_j!}\left(\frac{\sigma(j)+I_4(j)\sigma(j/4)}{j}\right)^{k_j}\right]\\
\nonumber&+&\sum_{m=1}^{n}\sum_{\pi(m)}(-5)^{\sum
k_j}\prod_{j=1}^{m}\left[\frac{1}{k_j!}\left(\frac{I_2(j)\sigma(j/2)}{j}\right)^{k_j}\right]\\
&\times&\sum_{\pi(n-m)}2^{\sum
k_j}\prod_{j=1}^{n-m}\left[\frac{1}{k_j!}\left(\frac{\sigma(j)+I_4(j)\sigma(j/4)}{j}\right)^{k_j}\right].
\end{eqnarray}

\end{document}